\documentclass[12pt,twoside]{article}
\usepackage{graphicx} % Заменен пакет graphics на graphicx
\usepackage{amsfonts,amssymb,amsmath} % подключение ресурсов ams-tex
\usepackage{euscript} % фигурные буквы

\date{}

\begin{document}

\centerline{}

\centerline{}

\centerline{ {\bf ON A FORMULATION OF THE CENTRAL LIMIT THEOREM}}

\centerline{}

\centerline{\bf {Alexander Shmyrov and Vasily Shmyrov}}

\centerline{}

\centerline{St. Petersburg State University}

\newtheorem{Theorem}{\quad Theorem}[section]
\newtheorem{Definition}[Theorem]{\quad Definition}
\newtheorem{Corollary}[Theorem]{\quad Corollary}
\newtheorem{Lemma}[Theorem]{\quad Lemma}
\newtheorem{Example}[Theorem]{\quad Example}

\centerline{}

\begin{abstract}
A non-classical formulation of the central limit theorem is given for sequences of independent random variables with finite second moments.
Singular sequences whose members all have a degenerate or normal distribution are excluded from consideration.
The condition of uniform convergence is imposed on the improper integrals defining the variances.
Under the conditions of non-singularity and uniform convergence, the central limit theorem is valid if and only if the total variance increases indefinitely.
\end{abstract}


\bigskip

\noindent {\bf Introduction}

The central limit theorem is a series of results on the convergence over the distribution of normalized sums of sequences of independent random variables with finite second moments to the normal standard distribution. The history of the central limit theorem is the history of the classes of such
sequences, which are determined by different conditions. These conditions and their combinations are imposed on the convergence.

Initially, the meaning of these conditions (Lyapunov, Lindeberg, and others) was to describe
the widest possible family of sequences of random variables for which it is guaranteed
convergence in distribution of sequences of normalized sums to a normal standard
random variable. It soon became clear that in order to obtain the necessary and sufficient conditions for such convergence, it is reasonable restrict ourselves to a certain class of sequences. For example, these are sequences satisfying the Kolmogorov (uniform neglect) condition [].

On the other hand, A.M. Lyapunov in his famous work "A new form of the theorem on the limit of probabilities" posed the description problem  of a class of sequences of random variables for which the central limit theorem becomes valid as soon as the total variance increases infinitely.
In this paper, it will be shown that such class is formed using a rather natural restriction, which makes it especially simple to formulate the central limit theorem.

\medskip
\noindent {\bf 1.} Let $\{\xi_n\}$ be a sequence of independent random variables with finite second
moments. We denote, as usual, by $ P_n $ - the probability distribution of a random variable
$ \ xi_n $, $ F_n $ - distribution function.

\noindent
Denote by:

$$
m_n \mathop = \limits^\Delta E \xi_n, \;\; \sigma_n^2 \mathop = \limits^\Delta D \xi_n, \;\; n = 1, \ldots,
$$

-- mathematical expectation and variance of a random variable $\xi_n$,
$$
S_n \mathop = \limits^\Delta \xi_1 + \cdots + \xi_n
$$
-- the sum of the first $n$ members of the sequence

$$
B_n^2 \mathop = \limits^\Delta D S_n = \sigma_1^2 + \cdots + \sigma_n^2
$$
-- total variance,

$$
\widetilde S_n  \mathop = \limits^\Delta \frac{S_n  - E S_n}{B_n}
$$
-- normalized sum.

For definiteness, we assume that the random variables $\xi_n$ are given on a countable line
product of sample spaces $(R,  \EuScript B(R), P_n)$ with Kolmogorov probability measure

$$
{\bf P} = P_1 \times \cdots \times P_n \times \cdots \, .
$$

Sequences of independent random variables with finite second moments are classical objects \cite{Lyapunov} of research in connection with different formulations of the central limit theorems.

The central limit theorem is a series of results about
convergence in distribution of sequences of normalized sums to a normal standard
random variable

$$
\tilde S_n  \mathop \to^d  \EuScript N (0,1).
$$

An interesting enough problem is to construct a class of sequences of random variables for which the central limit theorem becomes valid as soon as the total variance
increases infinitely.

In this paper, we show that such class is formed using a very natural restriction, which makes it possible to give a transparent formulation of the central
limit theorem \cite{Shmyrov}.

\medskip

\noindent {\bf 2}. The condition of finiteness of the second moments of the random variables $\xi_n$ means the convergence of improper integrals of the form

\begin{equation}
\label{Eq1}
\sigma_n^2 = \int  (x-m_n)^2 d F_n (x), \;\;\;  n = 1, \ldots,
\end{equation}
and one of the most common characteristics of sequences of the form \eqref{Eq1} in the analysis
is related to the concept of uniform convergence.
We use this concept in the case under consideration.

With the sequence \eqref{Eq1} we associate the sequence of functions $\{ \alpha_n \} $ of nonnegative argument $s$, which we define by the form

\begin{equation}
\label{Eq2}
\begin{array}{l}
\alpha_n (s) = 0, \;\; if \; \sigma_n =0; \\
\\
\alpha_n (s) = \frac{1}{\sigma_n^2} \int \limits_{|x-m_n| \ge s}  |x-m_n|^2  d F_n (s), \;\; if \; \sigma_n >0.
\end{array}
\end{equation}

It is clear, if the random variable $\xi_n$ has a finite second moment, then

$$
\alpha_n \to 0, \;\; s \to \infty.
$$

\medskip
\noindent {\bf Definition 1.} {\it If the random variables have finite second moments, and the functions $\{\alpha_n (s) \} $,defined by the formula \eqref{Eq2}, as $ s \to \infty $ tend to zero uniformly in $n$, then for a sequence of random variables there is the {\it condition of uniform convergence}. i.e.

\begin{equation}
\label{Eq3}
\mathop {\sup} \limits_n  \alpha_n (s) \to 0, \;\;\; s \to \infty.
\end{equation}

\smallskip
}
\medskip

The uniform convergence condition is a natural and not too restrictive complement to
the condition of finiteness of the second moments.

Obviously, the condition of uniform convergence is satisfied, if the set of distributions $\{ P_n \}$ of random variables $\{ \xi_n \}$ consists of a finite number
elements (with finite second moments, of course), as well as for uniformly bounded
random variables.

\bigskip

\noindent {\bf 3}. We formulate the central limit theorem using the uniform convergence condition for non-singular sequences.

\medskip

\noindent {\bf Definition 2.} {\it A sequence of random variables is called singular if all its members have a degenerate or normal distribution.}
\medskip

It is clear that if the sequence $\{\xi_n\}$ is singular, then the sequence of sums $\{S_n\}$ is also singular, so only the non-singular case remains interesting.

\medskip
\noindent
{\bf Theorem 1. The central limit theorem.}
{\it Let a non-singular sequence of random variables with finite second moments satisfy the condition of uniform convergence.
Then, in order for the sequence of normalized sums to converge in distribution to a normal standard random variable, it is necessary and sufficient that the total variance increases indefinitely.}

\medskip
\noindent {\bf P r o o f \, o f \, N e c e s s i t y.}
Without loss of generality, one can
assume that the random variable $\xi_1$ does not have a normal or degenerate distribution and that
the mathematical expectations of the random variables $\xi_n$ are equal to zero.

Let the total variance $ B_n^2 $ tends to a finite limit if $ n \to \infty $. Without loss of generality, this finite limit, , is set equal to one, i.e.

$$
B_n \to 1, \;\;\;   n \to \infty.
$$
We will show that then

$$
\widetilde S_n  \mathop {\not \to} \limits^d  \EuScript N (0,1).
$$
We have

$$
\widetilde S_n = \frac{S_n}{B_n} = \frac{\xi_1}{B_n} + \frac{\xi_2 + \cdots + \xi_n}{ B_n }
$$

and

$$
\frac{\xi_1}{B_n}  \mathop { \to} \limits^{ \bf P}  \xi_1.
$$

A sequence of random variables

$$
\eta_n \mathop = \limits^\Delta  \frac{\xi_2 + \cdots + \xi_n}{ B_n }, \;\;\; n= 1,2,\ldots,
$$
is fundamental in probability \cite{Shiryaev} since for any $\varepsilon> 0$ for $l>n$

$$
{\bf P}(| \eta_l -\eta_n |> \varepsilon) =
$$

$$
= {\bf P} \left (  \left |  \sum \limits_{k=1}^n \left (  \frac{1}{B_n} -\frac{1}{B_l}  \right ) \xi_k  +  \sum \limits_{k=n+1}^l  \frac{1}{B_l} \xi_k  \right | > \varepsilon  \right ) \leq
$$

$$
\leq   {\bf P} \left (  \left | \left (  \frac{1}{B_n} -\frac{1}{B_l}  \right )  \sum \limits_{k=1}^n   \xi_k   \right | > \frac{\varepsilon}{2} \right ) +
$$

$$
+ {\bf P} \left (  \left |  \frac{1}{B_l} \sum \limits_{k=n+1}^l \xi_k \right | > \frac{\varepsilon}{2} \right ) \leq  \left (  \frac{1}{B_l} -\frac{1}{B_n}  \right )^2  \frac{4 B_n^2}{\varepsilon^2}  +
$$

$$
+ \frac{ 4(B_l^2 - B_n^2) }{B_l^2 \varepsilon^2}  \to 0, \;\;\; n \to \infty.
$$

Since the sequence $\{ \eta_n \}$ is fundamental in probability, there is \cite{Shiryaev}
a random variable $\eta$ such that

$$
\eta_n \mathop { \to} \limits^{\bf P} \eta,
$$
moreover, the random variables $\xi_1$ and $\eta$ are independent. Thus,

$$
\widetilde S_n \mathop { \to} \limits^{\bf P}  \xi_1 + \eta
$$
and especially
$$
\widetilde S_n \mathop { \to} \limits^d  \xi_1 + \eta.
$$

If now
$$
\widetilde S_n  \mathop {\to} \limits^d  \EuScript N (0,1),
$$
then the random variable $\xi_1 + \eta$ has a normal distribution and by Cramer's theorem the random variables $\xi_1$ and $\eta$ have a normal distribution also. This contradicts the original assumption.

\medskip
\noindent
{\bf  P r o o f \,  o f \, s u f f i c i e n c y. }
Let
$$
B_n \to \infty.
$$

\noindent
Then the condition of uniform convergence \eqref{Eq3} implies the Lindeberg condition \cite{Lindeberg}, \cite{Lamperti}, since

$$
\frac{1}{B_n^2} \sum \limits_{k=1}^n  \int \limits_{|x-m_k| \ge \varepsilon B_n}  (x-m_k)^2 d F_k (x)  =
$$

$$
= \frac{1}{B_n^2} \sum \limits_{k=1}^n  \sigma_k^2  \alpha_k (\varepsilon B_n)  \le
$$

$$
\le  \mathop {\sup} \limits_k  \alpha_k (\varepsilon B_n)  \to 0, \;\;\; n \to \infty.
$$

It is known that Lindeberg condition guarantees that the sequence of distribution functions of normalized sums $ \{F_{\widetilde S_n} (x) \}$ converges to a normal standard function distribution uniformly on $x$. $\blacksquare$

\bigskip 
\noindent
{\bf e-mail: v.shmyrov@spbu.ru, vasilyshmyrov@yandex.ru}

\end{document}